\documentclass{article}

\usepackage[numbers,square]{natbib}
\usepackage{graphicx}
\usepackage{amsmath}
\usepackage{amsthm}
\usepackage{amssymb}

\newcommand{\E}{\mathbb E}
\newcommand{\R}{\mathbb{R}}
\newcommand{\N}{\mathbb{N}}

\renewcommand{\P}{\mathbb{P}}

\newcommand{\Var}{\mathop{\mathrm{Var}}\nolimits}

\newcommand{\eps}{\varepsilon}

\theoremstyle{plain}
\newtheorem{theorem}{Theorem}[section]

\theoremstyle{definition}

\theoremstyle{remark}
\newtheorem{remark}{Remark}[section]
\newtheorem{example}{Example}[section]

\begin{document}
\title{The maximum of the Gaussian $1/f^{\alpha}$-noise in the case $\alpha<1$}
\author{Zakhar Kabluchko}
\maketitle
\begin{abstract}
We prove that the appropriately normalized maximum of the Gaussian $1/f^{\alpha}$-noise with $\alpha<1$ converges in distribution to the Gumbel double-exponential law.
\end{abstract}

\section{Introduction and statement of the result}\label{sec:main}
$1/f^{\alpha}$-noise is usually described as a stochastic process whose spectral density is inverse proportional to some power of the frequency.
$1/f^{\alpha}$-noises have been observed experimentally in a huge variety of physical, biological, economic systems and are believed to be ubiquitous in nature. We refer  to~\cite{li_bib} for a list of references. The aim of the present paper is to find the limiting distribution of the maximum of the $1/f^{\alpha}$-noise in the case $\alpha<1$.

Let us be more precise. We define  $1/f^{\alpha}$-noise to be a Gaussian process $\{X_n(t), t\in [-\pi,\pi]\}$ given by a finite random Fourier series
\begin{equation}\label{eq:def_X_n}
X_n(t)=\sum_{k=1}^n \sqrt{R(k)} (U_k\sin (k t) +V_k \cos (k t)),
\end{equation}
where $U_k,V_k$ are independent real-valued standard Gaussian random variables and $R$ is some function  regularly varying at $+\infty$ with index $-\alpha$. We will recall necessary facts about regularly varying functions in Section~\ref{sec:reg_var}, the main example to keep in mind being $R(t)=ct^{-\alpha}$, where $c>0$.
Here we will be interested in the case $\alpha<1$.
In this case, for every $t\in[-\pi,\pi]$ the series on the right-hand side of~\eqref{eq:def_X_n} diverges as $n\to\infty$ with probability $1$. The next theorem is our main result.


\begin{theorem}\label{theo:main}
Let $X_n$ be the $1/f^{\alpha}$-noise defined by~\eqref{eq:def_X_n}, where $R:(0,\infty)\to[0,\infty)$ is an eventually monotone,  regularly varying function with index
$-\alpha$, where $-\infty<\alpha<1$. Let $\sigma_n^2=\Var X_n(0)$ and $c=2\pi^2\frac{1-\alpha}{3-\alpha}$. Then, for every $z\in\R$,
\begin{equation*}
\lim_{n\to\infty}\P\left[\frac 1{\sigma_n}\sup_{t\in[-\pi,\pi]}X_n(t)\leq
\sqrt {2\log n}+\frac{1}{\sqrt{2\log
n}}\left(\log \frac {\sqrt c}{\sqrt 2 \pi} + z\right) \right]=e^{-e^{-z}}.
\end{equation*}
\end{theorem}
\begin{example}
Taking $R(t)=1$, we obtain a limit theorem for the maximum of a random trigonometric polynomial $X_n(t)=\sum_{k=1}^n(U_k\sin (k t) +V_k \cos (k t))$.
\end{example}
\begin{remark}
The assumption $\alpha<1$ is crucial for the validity of Theorem~\ref{theo:main}. The case $\alpha>1$ is not interesting since in this case the series on the right-hand side of~\eqref{eq:def_X_n} converges uniformly with probability $1$; see~\cite[Ch.~VII,~\S 1,2]{kahane_book}.
This immediately implies that the maximum of $X_n$ converges weakly
(without any normalization) to the maximum of the corresponding infinite series. Much more interesting is the case $\alpha=1$. A non-Gumbel limiting distribution for the maximum of the $1/f$-noise has been derived by non-rigorous methods in the physical literature~\cite{fyodorov}. The maximum of the $1/f$-noise is believed to behave similarly to the maxima of other ``logarithmically correlated'' fields including the two-dimensional discrete Gaussian Free Field and the Branching Brownian Motion.  It has been shown recently that the maximum of the two-dimensional Gaussian Free Field recentered by its mean is tight~\cite{bramson}. It seems that the methods of~\cite{bramson} can be applied to the $1/f$-noise, but we will not do this here.
\end{remark}
The rest of the paper is devoted to the proof of Theorem~\ref{theo:main}. Throughout, $C$ is a large positive constant whose value may change from line to line.


\section{Method of the proof}\label{sec:method}

The idea of our proof of Theorem~\ref{theo:main} is to rescale the $1/f^{\alpha}$-noise in time in such a way that it becomes close to a stationary Gaussian process with differentiable sample paths. The limiting distribution for the maximum of such processes is recalled in the next theorem, see~\cite[Thm. 8.2.7]{leadbetter_book}.
\begin{theorem}[\cite{leadbetter_book}]\label{theo:lead}
Let $\{\xi(t), t\in\R\}$ be a stationary zero-mean, unit-variance Gaussian process with a.s.\ continuous paths. Suppose that the covariance function $\rho(t)=\E[\xi(0)\xi(t)]$ satisfies the following three conditions:
\begin{enumerate}
\item \label{a:1} For some $c>0$, $\rho(t)=1-ct^2+o(t^2)$ as $t\to 0$.
\item \label{a:2} $\lim_{t\to\infty} \rho(t)\log t=0$.
\item \label{a:3} $\rho(t)< 1$ for $t\neq 0$.
\end{enumerate}
Then, for every $z\in\R$,
$$
\lim_{n\to\infty}\P\left[\sup_{t\in [0,n]}\xi(t)\leq \sqrt {2\log n}+\frac{1}{\sqrt{2\log
n}}\left(\log \frac {\sqrt c}{\sqrt 2 \pi} + z\right) \right]=e^{-e^{-z}}.
$$
\end{theorem}

The following generalization of the above result to \textit{sequences} of
stationary Gaussian processes is due to~\citet{seleznjev}.
\begin{theorem}[\cite{seleznjev}]\label{theo:selez}
For every $n\in\N$ let $\{\xi_n(t),t\in [-\frac n2, \frac n2]\}$ be a stationary zero-mean, unit-variance
Gaussian process with a.s.\ continuous paths and covariance function $\rho_n(t)=\E[\xi_n(0)\xi_n(t)]$. Suppose that
\begin{enumerate}
\item  \label{c:1} $\rho_n(t)=1-c_n t^{2}+\eps_n(t), $
where $c_n$ is a sequence satisfying $\lim_{n\to\infty}c_n=c>0$  and $\eps_n(t)$ is a sequence of functions  satisfying $\lim_{t\to 0}\eps_n(t)/t^2=0$ uniformly in $n\in\N$.
\item \label{c:2} For every $\eps>0$ there is $T=T(\eps)$ such that
$\rho_n(t)\log t<\eps$ for every $n\in\N$, $t\in[T(\eps),\frac n2]$.
\item \label{c:3} For some $n_0\in\N$ and every $\eps>0$ we have $\sup_{n>n_0, t\in[\eps,\frac n2]}\rho_n(t)<1$.
\end{enumerate}
Then, for every $z\in\R$,
$$
\lim_{n\to\infty} \P\left[\sup_{t\in[-\frac n2,\frac n2]} \xi_n(t)\leq \sqrt {2\log n}+\frac{1}{\sqrt{2\log
n}}\left(\log \frac {\sqrt c}{\sqrt 2 \pi} + z\right)
\right] =e^{-e^{-z}}.
$$
\end{theorem}
Note that the conditions of Theorem~\ref{theo:selez} are just uniform versions of the conditions of Theorem~\ref{theo:lead}. An application of Theorem~\ref{theo:selez} can be found in~\cite{kabluchko}.

\section{Facts about regularly varying functions}\label{sec:reg_var}
We need to recall some facts from the theory of regular variation; see~\cite{bingham_book}. A positive measurable function $f$ defined on the positive half-axis is called \textit{regularly varying} at $+\infty$ with index $\alpha\in\R$ (notation: $f\in \text{RV}_{\alpha}$) if for every $\lambda>0$,
\begin{equation}\label{eq:def_reg_var}
\lim_{x\to+\infty}\frac{f(\lambda x)}{f(x)}=\lambda^{\alpha}.
\end{equation}
For example, the function $f(t)=ct^{\alpha}$, where $c>0$, is regularly varying with index $\alpha$. A regularly varying function with index $\alpha=0$ is called \textit{slowly varying}. Any function $f\in \text{RV}_{\alpha}$ can be written in the form $f(t)=L(t) t^{\alpha}$, where $L$ is slowly varying.

We will several times need the following result of Karamata~\cite[Prop.~1.5.8]{bingham_book}: if $f\in \text{RV}_{\alpha}$ with $\alpha>-1$, then
\begin{equation}\label{eq:karamata_sum}
\sum_{k=1}^n f(k) \sim \frac{nf(n)}{(1+\alpha)}, \;\;\; n\to\infty.
\end{equation}
(Note that Karamata's theorem is usually stated for the integral $\int_1^n f(t)dt$, but the discrete version given above is also true).
Also, wee will need an estimate called Potter bound~\cite[Thm.~1.5.6]{bingham_book}: if $L$ is slowly varying and bounded away from $0$ and $\infty$ on every compact subset, then for every $\delta>0$ there is a $C>0$ such that
\begin{equation}\label{eq:potter_bound}
\frac{L(x)}{L(y)}\leq C\max \left(\left(\frac xy\right)^{\delta}, \left(\frac yx\right)^{\delta}\right), \;\;\; x,y>0.
\end{equation}

\section{Proof of the main result}\label{sec:proof}
Let $X_n$ be a $1/f^{\alpha}$-noise as in Theorem~\ref{theo:main}. We represent the regularly varying function $R$ in the form $R(t)=L(t) t^{-\alpha}$, where $L$ is slowly varying.  The covariance function $r_n(t,s)=\E[X_n(t)X_n(s)]$ of the process $X_n$ is given by
\begin{equation}
r_n(t,s)=\sum_{k=1}^n R(k)\cos (k(t-s)), \;\;\;t,s\in [-\pi,\pi].
\end{equation}
In particular, for the variance $\sigma^2_n=\Var X_n(0)$ we have
\begin{equation}\label{eq:asympt_sigma_n}
\sigma_n^2=\sum_{k=1}^n R(k) \sim \frac{nR(n)}{(1-\alpha)}, \;\;\; n\to\infty,
\end{equation}
where the last step is a consequence of~\eqref{eq:karamata_sum} and the assumption $R\in \text{RV}_{-\alpha}$ with $\alpha<1$.
For every $n\in\N$ consider a rescaled process $\xi_n$ defined by
\begin{equation}\label{def:xi_n}
\xi_n(t)=\frac 1{\sigma_n} X_n\left(\frac {2\pi t}n\right),\;\;\; t\in \left[-\frac n2,\frac n2\right].
\end{equation}
Note that $\xi_n$ is a stationary Gaussian process with zero-mean, unit-variance margins. Its
 covariance function  $\rho_n(t)=\E [\xi_n(0)\xi_n(t)]$ is given by
\begin{equation}\label{eq:def_rho}
\rho_n(t)
=\frac 1{\sigma_n^2} \sum_{k=1}^n R(k)\cos \frac{2\pi kt}n, \;\;\; t\in \left[-\frac n2,\frac n2\right].
\end{equation}
We claim that the sequence $\xi_n$, $n\in\N$, satisfies the assumptions of Theorem~\ref{theo:selez}. We start by verifying condition~\ref{c:1}.
Write
\begin{equation}
\delta_{k,n}(t)=\cos\frac{2\pi kt}{n}-\left(1-\frac{2\pi^2 k^2t^2}{n^2}\right).
\end{equation}
Then,
\begin{equation}\label{eq:rho_rep}
\rho_n(t)
=
1-\frac{2\pi^2 t^2}{n^2\sigma_n^{2}} \sum_{k=1}^n k^2R(k)+\frac 1{\sigma_n^{2}}\sum_{k=1}^n R(k)\delta_{k,n}(t).
\end{equation}
Note that the function $t^2R(t)$ is regularly varying with index $2-\alpha>-1$. By~\eqref{eq:karamata_sum} and~\eqref{eq:asympt_sigma_n} we have
\begin{equation}\label{eq:c_n}
c_n:= \frac{2\pi^2}{n^2\sigma_n^{2}} \sum_{k=1}^n k^2R(k)\to 2\pi^2\frac{1-\alpha}{3-\alpha}, \;\;\; n\to\infty.
\end{equation}
Let us estimate the third term on the right-hand side of~\eqref{eq:rho_rep}.
By Taylor's expansion, for every $\eps>0$ there exists a $\delta>0$ such that $|\delta_{k,n}(t)|<\eps t^2$ for every $n\in\N$, $1\leq k\leq n$ and $|t|<\delta$. It follows that
\begin{equation}
|\eps_n(t)|
:=
\left|\frac 1{\sigma_n^{2}}\sum_{k=1}^n R(k)\delta_{k,n}(t)\right|
\leq
\frac{\eps t^2}{\sigma_n^2} \sum_{k=1}^n R(k)
=\eps t^2
\end{equation}
uniformly over $n\in\N$,  $|t|<\delta$.
Together with~\eqref{eq:rho_rep} and~\eqref{eq:c_n} this proves that condition~\ref{c:1} of Theorem~\ref{theo:selez} holds with $c=2\pi^2\frac{1-\alpha}{3-\alpha}$.



Let us now show that condition~\ref{c:2} of Theorem~\ref{theo:selez} is satisfied. To estimate $\rho_n(t)$ for large $t$ we need to take into account the oscillating character of the terms on the right-hand side of~\eqref{eq:def_rho}, which suggests performing Abel's summation.  However, it can be shown that a direct application of Abel's summation  leads to a satisfactory estimate for $\alpha<0$ only.  So, we need a somewhat more accurate argument.

First of all, we may redefine the function $R$ on an interval of the form $(0,A)$ to make it monotone on the whole positive half-line and bounded away from $0$ on any compact set. Indeed, such a modification changes $\rho_n(t)$ by at most $C/\sigma_n^2$ (see~\eqref{eq:def_rho}) which is smaller than $\eps/(2\log t)$ uniformly in $t\in [T(\eps),\frac n2]$, where $T(\eps)$ is large. So, the modification has no influence on the validity of condition~\ref{c:2} of Theorem~\ref{theo:selez}.

Let $1\leq t\leq \frac n2$. We will split the sum defining $\rho_n(t)$ as follows:
\begin{equation}\label{eq:est1}
\rho_n(t)
=\frac 1{\sigma_n^2} \sum_{k=1}^{[n/t]-1} R(k)\cos \frac{2\pi kt}{n}+
\frac 1{\sigma_n^2} \sum_{k=[n/t]}^{n} R(k)\cos \frac{2\pi kt}{n}=:S_1+S_2.
\end{equation}
The sum $S_1$ can be estimated in a trivial way: using the inequality $|\cos x|\leq 1$ and~\eqref{eq:asympt_sigma_n}, we obtain
\begin{equation}\label{eq:est2}
|S_1|
\leq \frac 1{\sigma_n^2} \sum_{k=1}^{[n/t]} R(k)
=\frac{\sigma^2_{[n/t]}}{\sigma_n^2}
\leq C\frac{(n/t) R(n/t)}{nR(n)} = Ct^{\alpha-1}\frac{L(n/t)}{L(n)}.
\end{equation}
The sum $S_2$  will be estimated by Abel's summation. We need the Dirichlet kernel
\begin{equation}
D_k(t)=\sum_{j=1}^k \cos \left(2\pi jt\right)=\frac 12 \left(\frac{\sin((2k+1)\pi t)}{\sin(\pi t) }-1\right).
\end{equation}
Since $|t|\leq \frac n2$, we have $|\sin\frac{\pi t} n|\geq \kappa \frac tn$ for some $\kappa>0$. It follows that for every $k\in\N$,
\begin{equation}\label{eq:dir_kern_est}
\left|D_k\left(\frac tn\right)\right|\leq \frac 12+ \frac 1 {\kappa}\frac nt\leq C\frac nt.
\end{equation}
Applying Abel's summation formula to the sum $S_2$ we obtain
\begin{align}
S_2
&=
\frac 1{\sigma_n^{2}}\sum_{k=[n/t]}^{n-1} D_k\left(\frac tn\right) \left(R(k)-R(k+1) \right)\notag\\
&+\frac 1{\sigma_n^{2}} D_n\left(\frac tn\right) R(n)
-\frac 1{\sigma_n^{2}} D_{[n/t]-1}\left(\frac tn\right) R([n/t]).\label{eq:rho_abel_summ}
\end{align}
Utilizing \eqref{eq:dir_kern_est} and~\eqref{eq:asympt_sigma_n} we get that $|S_2|$ can be estimated from above by
\begin{align}
|S_2|&\leq
\frac 1{\sigma_n^{2}}\sum_{k=[n/t]}^{n-1} \left|D_k\left(\frac tn\right)\right|\left|R(k)-R(k+1) \right|\notag\\
&+\frac 1{\sigma_n^{2}} \left|D_n\left(\frac tn\right)\right|R(n) +\frac 1{\sigma_n^{2}} \left|D_{[n/t]-1}\left(\frac tn\right)\right| R([n/t]) \notag\\
&\leq
\frac{C}{tR(n)}\left(\sum_{k=[n/t]}^{n-1} \left|R(k)-R(k+1) \right|+ R(n)+R(n/t)\right).\label{eq:est_S2}
\end{align}
Recall that $R$ is assumed to be monotone. Depending on whether $R$ is decreasing or increasing, the expression in the brackets in~\eqref{eq:est_S2} can be estimated from above by $2R(n/t)$ or $2R(n)$. Thus,
\begin{equation}\label{eq:est3}
|S_2|\leq
C\max\left(\frac {R(n/t)}{tR(n)},\frac 1t\right)
=
C\max\left(t^{\alpha-1}\frac{L(n/t)}{L(n)},\frac 1t\right).
\end{equation}
Bringing~\eqref{eq:est1}, \eqref{eq:est2}, \eqref{eq:est3} together and employing Potter's bound~\eqref{eq:potter_bound} we obtain that for every $\delta>0$ there is $C>0$ such that for all $n\in\N$ and $t\in[1,\frac n2]$,
\begin{equation}
|\rho_n(t)|
\leq C\max\left(t^{\alpha-1}\frac{L(n/t)}{L(n)},\frac 1t\right)
\leq C\max\left( t^{\alpha-1+\delta}, \frac 1t\right).
\end{equation}
Recall that we assume that $\alpha<1$. Choose $\delta>0$ so small that $\alpha-1+\delta<0$. The verification of condition~\ref{c:2} is completed.


Let us finally verify condition~\ref{c:3} of Theorem~\ref{theo:selez}. Fix $\eps>0$. By condition~\ref{c:2} there is $T>0$ such that for all $n\in\N$, $t\in [T, \frac n2]$ we have $\rho_n(t)<1/2$. Thus, we have to show that for some $n_0\in\N$,
\begin{equation}\label{eq:need_cond3}
\sup_{n>n_0, t\in [\eps,T]}\rho_n(t)<1.
\end{equation}
We can find sufficiently small $a>0$ and $\eta>0$ such that $\left|\cos \frac{2\pi kt}{n}\right|<1-\eta$ for all $n\in\N$, $k\in [an,2an]$, $t\in [\eps,T]$.
Recalling~\eqref{eq:def_rho} we have
\begin{equation}
\rho_n(t)
=
\frac 1{\sigma_n^2} \sum_{k\in [an, 2an]} R(k)\cos \frac{2\pi kt}{n}+
\frac 1{\sigma_n^2} \sum_{\substack{1\leq k\leq n\\k\notin [an, 2an]}} R(k)\cos \frac{2\pi kt}{n}.
\end{equation}
It follows that for all $n\in\N$ and $t\in [\eps,T]$,
\begin{equation}
|\rho_n(t)|\leq
\frac 1{\sigma_n^2} \sum_{k\in [an, 2an]} (1-\eta) R(k)
+
\frac 1{\sigma_n^2} \sum_{\substack{1\leq k\leq n\\k\notin [an, 2an]}} R(k)
=
1- \frac {\eta}{\sigma_n^2} \sum_{k\in [an, 2an]} R(k).
\end{equation}
Applying~\eqref{eq:karamata_sum} and~\eqref{eq:asympt_sigma_n}, we obtain that uniformly in $t\in [\eps,T]$,
\begin{equation}
\limsup_{n\to\infty}|\rho_n(t)|\leq  1- \eta \lim_{n\to\infty}\frac{2anR(2an)-anR(an)}{nR(n)}=1- (2^{1-\alpha}-1)a^{1-\alpha}\eta<1.
\end{equation}
Hence, there is $n_0\in\N$ such that~\eqref{eq:need_cond3} holds. This verifies condition~\ref{c:3} of Theorem~\ref{theo:selez}.

To complete the proof of Theorem~\ref{theo:main}, note that $ \sigma_n^{-1}\sup_{t\in[-\pi,\pi]}X_n(t)$ has the same law as $\sup_{t\in[-\frac n2,\frac n2]} \xi_n(t)$ and apply Theorem~\ref{theo:selez}.
\bibliographystyle{plainnat}
\bibliography{paper22bib}
\end{document}